\documentclass[12pt]{article}
\usepackage{e-jc}
\usepackage{amsmath}
\usepackage{amssymb, amsthm}
\newtheorem{theorem}{Theorem}[section]

\newtheorem{conj}[theorem]{Conjecture}

\newcommand{\erd}{Erd\H{o}s}
\newcommand{\multichoose}[2]{\left(\! \binom{#1}{#2}\! \right)}

\title{An Erd\H{o}s-Ko-Rado theorem for multisets}
\author{Karen Meagher \footnote{Research supported by NSERC.}\\
\small  Department of Mathematics and Statistics \\[-0.8ex]
\small University of Regina,  Regina, Saskatchewan, Canada\\[-0.8ex]
\small \texttt{kmeagher@math.uregina.ca}\\
Alison Purdy \\ [-0.8ex]
\small  Department of Mathematics and Statistics \\[-0.8ex]
\small University of Regina,  Regina, Saskatchewan, Canada\\[-0.8ex]
\small \texttt{purdyali@math.uregina.ca}
}

 \date{\dateline{Aug 9, 2011}{Oct 26, 2011}\\
   \small Mathematics Subject Classification: 05D05}
   
\begin{document}
\maketitle

\abstract{Let $k$ and $m$ be positive integers.  A collection of $k$-multisets from $\{1,\dots, m \}$ is intersecting if every pair of multisets from the collection is intersecting.  We prove that for $m \geq k+1$,  the size of the largest such collection is $\binom{m+k-2}{k-1}$ and that when $m > k+1$, only a collection of all the $k$-multisets containing a fixed element will attain this bound.  The size and structure of the largest intersecting collection of $k$-multisets for $m \leq k$ is also given.
\section{Introduction}

The \erd-Ko-Rado Theorem~\cite{MR0140419} is an important result in extremal set theory that gives the size and structure of the largest pairwise intersecting $k$-subset system from $[n] = \{1,\dots,n\}$.  This theorem is commonly stated as follows:  
\begin{theorem}
Let $k$ and $n$ be positive integers with $n \geq 2k$.  If $\mathcal F$ is a collection of intersecting $k$-subsets of $[n]$, then \[|\mathcal F| \leq \binom{n-1}{k-1}.\]  Moreover, if $n >2k$, equality holds if and only if $\mathcal F$ is a collection of all the $k$-subsets from $[n]$ that contain a fixed element from $[n]$.
\end{theorem}

Note that if $n < 2k$, any pair of $k$-subsets will be intersecting and so the largest intersecting collection will have size $\binom{n}{k}$.

A multiset is a generalization of a set in which an element may appear more than once.  As with sets, the order of the elements is irrelevant.  The cardinality of a multiset is the total number of elements including repetitions.  A $k$-multiset system on an $m$-set is a collection of multisets of cardinality $k$ containing elements from $[m]$.  We say that two multisets are intersecting if they have at least one element in common and that a collection of multisets is intersecting if every pair of multisets in the collection is intersecting. 

In this paper, we give a generalization of the \erd-Ko-Rado Theorem to intersecting multiset systems.  Specifically, we prove the following two theorems for the cases when $m \geq k+1$ and $m \leq k$ respectively.

\begin{theorem}\label{m>k}
Let $k$ and $m$ be positive integers with $m \geq k +1 $.  If $\mathcal A$ is a collection of intersecting $k$-multisets of $[m]$, then \[|\mathcal A| \leq \binom{m + k - 2}{ k-1}.\]  Moreover, if $m >k + 1$, equality holds if and only if $\mathcal A$ is a collection of all the $k$-multisets from $[m]$ that contain a fixed element from $[m]$.
\end{theorem}

If $m< k+1$, larger collections are possible.  For example, if $m\!=\!k\!=\!3$, the seven $k$-multisets containing either two or three distinct elements from $[m]$ will form an intersecting collection since each multiset contains more than half the elements from $[m]$.  We will use $\mathcal M_{(> \frac{m}{2})}$ to denote the collection of all $k$-multisets that contain more than $\frac{m}{2} $ distinct elements from $[m]$ and $\mathcal M_{(\frac{m}{2})}$ to denote the collection of all $k$-multisets from $[m]$ containing exactly $\frac {m}{2}$ distinct elements.  Then $$\left| \mathcal M_{(\frac{m}{2})} \right|=\binom{m}{\frac{m}{2}}\binom{k-1}{k-\frac{m}{2}}$$ and $$\left|\mathcal M_{(> \frac{m}{2})}\right| = \sum_{j=\lceil \frac{m+1}{2}\rceil}^{m}\binom{m}{j}\binom{k-1}{k-j}.$$

\begin{theorem}\label{m<k+1}
Let $k$ and $m$ be positive integers with $m \leq k$.  If $\mathcal A$ is a collection of intersecting $k$-multisets of $[m]$, then:
\begin{enumerate}
\item If $m$ is odd, $\left|\mathcal A \right| \leq \left|\mathcal M_{(> \frac{m}{2})} \right|$
and equality holds if and only if $\mathcal A = \mathcal M_{(> \frac{m}{2})}$. 
\item  If $m$ is even, $\left |\mathcal A \right| \leq \left|\mathcal M_{(> \frac{m}{2})}\right| + \frac{1}{2} \left | \mathcal M_{(\frac{m}{2})} \right|$
and equality holds if and only if $\mathcal A$ consists of $\mathcal M_{(>\frac{m}{2})}$ and a maximal intersecting collection of  $k$-multisets from $\mathcal M_{(\frac{m}{2})}$.
\end{enumerate}
\end{theorem}

A $k$-multiset on an $m$-set can be represented as an integer sequence of length $m$ with the integer in each position representing the number of repetitions of the corresponding element from $[m]$.  For example, if $m=6$, the multiset $\{1,2,2,4\}$ can be represented by the integer sequence $(1,2,0,1,0,0)$.  For a $k$-multiset, the sum of the integers in the corresponding integer sequence will equal $k$.

\erd-Ko-Rado type results for intersecting families of integer sequences are known (e.g.~\cite{MR593848},  \cite{MR1722210}, \cite{MR721369}).  In these, the sum of the entries in the integer sequence is not restricted to $k$ and the definition of intersection is different from our definition for multisets.  In~\cite{MR935022}, Anderson proves an \erd-Ko-Rado type result for multisets but uses yet another definition of intersection.  A definition of intersection equivalent to ours is used in several theorems for intersecting collections of vectors presented by Anderson in \cite{MR1902962}.  These theorems were originally written in terms of sets of noncoprime divisors of a number by Erd\H{o}s et al. in \cite{MR0274413} and \cite{MR0313065}, and again the sum of the entries is not restricted to $k$.  

More recently, Brockman and Kay~\cite{bk} proved the result in Theorem~\ref{m>k} provided that $m \geq 2k$.  Mahdian~\cite{mm} proved the bound on the size of a collection of intersecting $k$-multisets when $m >k$ using a method similar to Katona's cycle proof for sets~\cite{MR0304181}.  Our results improve the bound on $m$ given in~\cite{bk} and give the size and structure of the largest possible intersecting collection for all values of $m$ and $k$.

\section{Proof of Theorem~\ref{m>k}}

Our proof of this theorem uses a homomorphism from a Kneser graph to a graph whose vertices are the $k$-multisets of $[m]$.

A Kneser graph, $K(n,k)$, is a graph whose vertices are all of the $k$-sets from an $n$-set, denoted by $\binom{[n]}{k}$, and where two vertices are adjacent if and only if the corresponding $k$-sets are disjoint.  Thus an independent set of vertices in the Kneser graph is an intersecting $k$-set system.  We will use $\alpha (K(n,k))$ to denote the size of the largest independent set in $K(n,k)$.

We now define a multiset analogue of the Kneser graph.  For positive integers $k$ and $m$, let $M(m,k)$ be the graph whose vertices are the $k$-multisets from the set $[m]$, denoted by $\multichoose {[m]}{k}$, and where two vertices are adjacent if and only if the corresponding multisets are disjoint.  For this graph, the number of vertices is equal to $\multichoose {m}{k} = \binom{m+k-1}{k}$ and an independent set is an intersecting $k$-multiset system.

Let $n=m+k-1$.  Then $K(n,k)$ has the same number of vertices as $M(m,k)$ and $B \cap [m] \neq \emptyset$ for any $B \in \binom{[n]}{k}$. 

For a set $A \subseteq [m]$ of cardinality $a$ where $1 \leq a \leq k$, the number of $k$-sets, $B$, from $[n]$ such that $B \cap [m] = A$ will be equal to $$\binom{n-m}{k-a} = \binom{k-1}{k-a}.$$  Similarly, the number of $k$-multisets from $[m]$ that contain all of the elements of $A$ and no others will be equal to $$\multichoose{a}{k-a} =\binom{a+(k-a)-1}{k-a } = \binom{k-1 }{k-a}.$$

Hence there exists a bijection, $f: \binom{[n]}{k} \rightarrow \multichoose {[m]}{k}$, such that for any $B \in \binom{[n] }{k}$, the set of distinct elements in $f(B)$ will be equal to $B \cap [m]$.

If $A,B \in \binom{[n]}{k}$ are two adjacent vertices in the Kneser graph, then $(A \cap [m]) \cap (B \cap [m]) = \emptyset$ and hence $f(A) \cap f(B) = \emptyset$.  Therefore $f(A)$ is adjacent to $f(B)$ if $A$ is adjacent to $B$ and so the bijection $f: \binom{[n]}{k} \rightarrow \multichoose {[m]}{k}$ is a graph homomorphism.  In fact, $K(n,k)$ is isomorphic to a spanning subgraph of $M(m,k)$.  
Thus $$\alpha(M(m,k)) \leq \alpha (K(n,k)).$$
From the \erd-Ko-Rado Theorem, we have that if $n \geq 2k$,  $$\alpha (K(n,k)) = \binom{n-1}{k-1}.$$  Thus, for  $m \geq k+1$, $$\alpha(M(m,k)) \leq \binom{n-1}{k-1} = \binom{m+k-2}{k-1}.$$

An intersecting collection of $k$-multisets from $[m]$ consisting of all $k$-multisets containing a fixed element from $[m]$ will have size $\binom{m + (k-1) -1}{k-1}=\binom{m-k-2}{k-1}$.  Therefore $$\alpha (M(m,k)) = \binom{m+k-2}{k-1}$$ which gives the upper bound on $\mathcal A$ in Theorem~\ref{m>k}.

To prove the uniqueness statement in the theorem, let $m > k+1$ and let $\mathcal A$ be an intersecting multiset system of size $\binom{m+k-2}{k-1}$.  With the homomorphism defined above, the pre-image of $\mathcal A$ will be an independent set in $K(n,k)$ of size $\binom{n-1}{k-1}$.  Since $m >k+1$ and $n=m+k-1$, it follows that $n>2k$ so, by the \erd-Ko-Rado theorem, $f^{-1}(\mathcal A)$ will be a collection of all the $k$-subsets of $[n]$ that contain a fixed element from $[n]$.  If the fixed element, $x$, is an element of $[m]$, then it follows from the definition of $f$ that every multiset in $\mathcal A$ will contain $x$.  Thus $\mathcal A$ will be a collection of all the $k$-multisets from $[m]$ that contain a fixed element from $[m]$ as required.  If $x \notin [m]$, then $f^{-1}(\mathcal A)$ will include the sets $A=\{1,m+1,\dots,n\}$ and $B = \{2,m+1,\dots, n\}$ since $m > k+1$ implies that $m > 2$.  But $f(A) \cap f(B)=\emptyset$ which contradicts our assumption that $\mathcal A$ is an intersecting collection of multisets.  Therefore, when $m > k+1$, if $\mathcal A$ is an intersecting collection of multisets of the maximum possible size, then $\mathcal A$ is the collections of all $k$-multisets containing a fixed element from $[m]$.
\qed

The case when $m=k+1$ is analogous to the case when $n=2k$ in the \erd-Ko-Rado theorem.  The size of the largest possible intersecting collection is equal to $\binom{m+k-2}{k-1}$ but collections attaining this bound are not limited to those having a common element in all $k$-multisets.

\section{Proof of Theorem~\ref{m<k+1}}
Although Theorem~\ref{m>k} is restricted to $m \geq k+1$, the inequality $\alpha(M(m,k)) \leq \alpha (K(n,k))$ still holds when $m \leq k$.  However, the resulting inequality $$\alpha(M(m,k)) \leq \binom{n}{k}=\binom {m+k-1}{k}$$ is not particularly useful since for $m >1$ this bound is not attainable.  Clearly, two multisets consisting of $k$ copies of different elements from $[m]$ will not intersect.

Before proceeding with our proof of Theorem~\ref{m<k+1}, we define the support of a multiset.  If $A$ is a $k$-multiset from $[m]$, the support of $A$, denoted by $S_{A}$, is the set of distinct integers from $[m]$ in $A$.  Note that two $k$-multisets, $A,B \in \multichoose{[m]}{k}$, will be intersecting if and only if $S_{A} \cap S_{B} \neq \emptyset$ and that each $S_{A}$ will have a unique complement, $\overline{S_{A}}=[m]\backslash S_{A}$, in $[m]$.  

Let $\mathcal A$ be an intersecting family of $k$-multisets of $[m]$ of maximum size and let $M \in \mathcal A$ be a $k$-multiset such that $|S_{M}| = \min \{\,\left|S_{A}\right| : A \in \mathcal A\}$. If $m=2$, it is easy to see that the theorem holds, so we will assume that $m > 2$. 

Suppose that $\left| S_{M} \right| < \frac{m}{2}$.  Let $\mathcal B_{1}= \{A \in \mathcal A : S_{A}=S_{M}\}$ and let $\mathcal B_{2}=\{B \in \multichoose {[m]}{k} : S_{B}= \overline{S_{M}}\}$.  Then $\mathcal B_{1} \subseteq \mathcal A$ and  $\mathcal B_{2} \cap \mathcal A = \emptyset$.
 
We will now show that ${\mathcal A}' := (\mathcal A \backslash \mathcal B_{1}) \cup \mathcal B_{2}$ is an intersecting family of $k$-multisets from $[m]$ that is larger than $\mathcal A$.  By construction, every multiset in $\mathcal A \backslash \mathcal B_{1}$ contains at least one element from $[m]\backslash S_{M}$, and $[m]\backslash S_{M} =S_{B}$ for all $B \in \mathcal B_{2}$.  Thus ${\mathcal A}'$ is an intersecting collection of $k$-multisets.

Let $\left| S_{M} \right| = i$.  Then  

$$\left| \mathcal B_{1} \right| =  \multichoose {i}{k-i} = \binom{k-1}{k-i}.$$ 

Since $\left| \overline{S_{M}} \right| = m-i$, it follows that
$$\left| \mathcal B_{2} \right| = \multichoose {m-i}{k-(m-i)} = \binom{k-1}{k-m+i}.$$

To show that $\left| {\mathcal A}' \right| >  \left| \mathcal A \right|$, it is sufficient to show that $$ \binom{k-1}{ k-m+i} > \binom{k-1}{k-i},$$ or equivalently, that $$(k-i)!(i-1)! > (k-m+i)!(m-i-1)!\, .$$

Since $i < \frac{m}{2}$ and $m \leq k$, we have that $k-i > k-\frac{m}{2} > k-m+i \geq 1$.  Therefore, 
\begin{align*}
(k-i)!(i-1)! &= (k-i)(k-i-1)\dots(k-m+i+1)(k-m+i)!(i-1)!\\
&\geq (m-i)(m-i-1)\dots(i+1)(k-m+i)!(i-1)!\\
&=\frac{m-i}{i}(m-i-1)!(k-m+i)!\\
&> (m-i-1)!(k-m+i)!
\end{align*}
as required.  Thus if $\mathcal A$ is of maximum size, it cannot contain a multiset with less than $\frac{m}{2}$ distinct elements from $[m]$.

It is easy to see that any $k$-multiset containing more than $\frac{m}{2}$ distinct elements from $[m]$ will intersect with any other such $k$-multiset.  This completes the proof of the theorem for the case when $m$ is odd.  When $m$ is even, it is necessary to consider the $k$-multisets which contain exactly $\frac{m}{2}$ distinct elements, that is, the $k$-multisets in $\mathcal M_{(\frac{m}{2})}$.  These multisets will intersect with any multiset containing more than $\frac{m}{2}$ distinct elements.  However, $\mathcal M_{(\frac{m}{2})}$ is not an intersecting collection.  For any $A \in \mathcal M_{(\frac{m}{2})}$, all of the $k$-multisets, $B$, where $S_{B}=\overline{S_{A}}$ will be in $\mathcal M_{(\frac{m}{2})}$ and will not intersect with $A$.  Since the size of a maximal intersecting collection of $\frac{m}{2}$-subsets of $[m]$ is $\frac{1}{2}\binom{m}{\frac{m}{2}}$ and each $\frac{m}{2}$-subset is the support for the same number of multisets in $\mathcal M_{(\frac{m}{2})}$, an intersecting collection of $k$-multisets will contain at most half of the $k$-multisets in $\mathcal M_{(\frac{m}{2})}$.
\qed

\section{Further work}

An obvious open problem is determining the size and structure of the largest collection of $t$-intersecting $k$-multisets, i.e. collections of multisets where the size of the intersection for every pair of multisets is at least $t$.  (We define the intersection of two multisets to be the multiset containing all elements common to both multisets with repetitions.)  The following conjecture is a version of Conjecture 5.1 from~\cite{bk}.

\begin{conj}\label{conj}
Let $k$, $m$ and $t$ be positive integers with $t \leq k$ and $m \geq t(k-t) + 2 $.  If $\mathcal A$ is a collection of intersecting $k$-multisets of $[m]$, then \[|\mathcal A| \leq \binom{m + k - t - 1}{k-t}.\]  Moreover, if $m >t(k-t) + 2$, equality holds if and only if $\mathcal A$ is a collection of all the $k$-multisets from $[m]$ that contain a fixed $t$-multiset from $[m]$.
\end{conj}

The lower bound on $m$ in this conjecture was obtained by substituting $m+k-1$ for $n$ in the corresponding bound for sets given by Frankl~\cite{MR519277} and Wilson~\cite{MR771733}.  The conjecture is supported by the fact that when $m > t(k-t)+2 $, a collection consisting of all $k$-multisets containing a fixed $t$-multiset is larger than a collection consisting of all $k$-multisets containing $t+1$ elements from a set of $t+2$ distinct elements of $[m]$ and that these two collections are equal in size when $m=t(k-t) + 2$.  Furthermore, when $m = t(k-t) + 1$, collections larger than $\binom{m + k - t - 1}{k-t}$ are possible.  For example, if $t=2$, $k=5$ and $m=7$, the cardinality of the collection of all $k$-multisets containing three or more elements from $\{1,2,3,4\}$ is $91$ while $\binom{m + k - t - 1}{k-t}=84$.

The existence of a graph homomorphism from the Kneser graph $K(n,k)$ to its multiset analogue $M(m,k)$ in the proof of Theorem~\ref{m>k} gave a simple and straight-forward way to show that the size of the largest independent set in $M(m,k)$ is no larger than the size of the largest independent set in $K(n,k)$.  These graphs can be generalized as follows:  let $K(n,k,t)$ be the graph whose vertices are the $k$-subsets of $[n]$ and where two vertices, $A,B$, are adjacent if $\left| A \cap B \right| < t$ and let $M(m,k,t)$ be the graph whose vertices are the $k$-multisets of $[m]$ and where two vertices, $C,D$ are adjacent if  $\left| C \cap D \right| < t$.

If a bijective homomorphism from $K(n,k,t)$ to $M(m,k,t)$ exists, it could be used to prove a bound not only on the maximum size of a $t$-intersecting collection as given in Conjecture~\ref{conj} but also on the maximum size when $k-t \leq m \leq t(k-t) + 2$ using the Complete \erd-Ko-Rado theorem of Ahlswede and Khachatrian~\cite{MR1429238}.  However, it is not clear that such a homomorphism exists.  The conditions placed on the bijection in the proof of Theorem~\ref{m>k} are not sufficient to ensure that the bijection is a homomorphism since for two $k$-multisets, $A$ and $B$, having $\left| S_{A} \cap S_{B} \right| < t$ does not imply that $\left| A \cap B \right| < t$.

The simple fact that if a graph $G$ is isomorphic to a spanning subgraph of a graph $H$, then $\alpha(H) \leq \alpha(G)$
may be useful in proving \erd-Ko-Rado theorems for different objects.  It would be interesting to determine if there are  combinatorial objects other than multisets which have this relationship to an object for which an \erd-Ko-Rado type result is known. 
 
\section*{Acknowledgment}
We are grateful for the helpful comments of the anonymous referee, particularly those concerning Theorem~\ref{m<k+1} which greatly simplified the proof.

\end{document}